\documentclass[12pt, hidelinks]{article}

\usepackage[utf8]{inputenc}
\usepackage{url,hyperref}
\usepackage{amsmath,amssymb}
\usepackage{fullpage}
\usepackage{xcolor}
\usepackage[small, bf]{caption}
\usepackage{subcaption}
\usepackage{graphicx}
\graphicspath{{figures/}}
\usepackage{algorithm, algpseudocode}

\usepackage[style=alphabetic,backend=bibtex]{biblatex}
\bibliography{refs}
\newcommand{\BEAS}{\begin{eqnarray*}}
    \newcommand{\EEAS}{\end{eqnarray*}}
    \newcommand{\BEA}{\begin{eqnarray}}
    \newcommand{\EEA}{\end{eqnarray}}
    \newcommand{\BEQ}{\begin{equation}}
    \newcommand{\EEQ}{\end{equation}}
    \newcommand{\BIT}{\begin{itemize}}
    \newcommand{\EIT}{\end{itemize}}
    \newcommand{\BNUM}{\begin{enumerate}}
    \newcommand{\ENUM}{\end{enumerate}}
    
    \newcommand{\eg}{{\it e.g.}}
    \newcommand{\ie}{{\it i.e.}}

    \newcommand{\ones}{\mathbf 1}
    
    \newcommand{\reals}{\mathbb{R}}

    \newcommand{\bmat}[1]{\begin{bmatrix}#1\end{bmatrix}}

    \newif\iftodos

\usepackage[USenglish]{babel}
\usepackage[T1]{fontenc}

\usepackage[]{jlcode}

\title{An Efficient Algorithm for Optimal Routing Through Constant Function Market Makers}

\author{
    Theo Diamandis\\
    {\small \texttt{tdiamand@mit.edu}}
    \and
    Max Resnick\\
    {\small \texttt{max@riskharbor.com}}
    \and
    Tarun Chitra\\
    {\small \texttt{tarun@gauntlet.network}}
    \and
    Guillermo Angeris \\
    {\small \texttt{gangeris@baincapital.com}}
}
\date{February 2023}

\newcommand{\arb}{\mathbf{arb}}

\let\phi\varphi
\let\eps\varepsilon

\begin{document}

%--------------------------------- Text ----------------------------------

\maketitle

\begin{abstract}
Constant function market makers (CFMMs) such as Uniswap have facilitated
trillions of dollars of digital asset trades and have billions of dollars of
liquidity.
One natural question is how to optimally route trades across a network of CFMMs
in order to ensure the largest possible utility (as specified by a user). We
present an efficient algorithm, based on a decomposition method, to solve the
problem of optimally executing an order across a network of decentralized
exchanges. The decomposition method, as a side effect, makes it simple to
incorporate more complicated CFMMs, or even include `aggregate CFMMs' (such as
Uniswap v3), into the routing problem.
Numerical results show significant performance improvements of this method,
tested on realistic networks of CFMMs, when compared against an off-the-shelf
commercial solver.
\end{abstract}

\section*{Introduction}
Decentralized Finance, or DeFi, has been one of the largest growth areas within
both financial technologies and cryptocurrencies since 2019. DeFi is made up of
a network of decentralized protocols that match buyers and sellers of digital
goods in a trustless manner. Within DeFi, some of the most popular applications are
decentralized exchanges (DEXs, for short) which allow users to permissionlessly
trade assets. While there are many types of DEXs, the most popular form of
exchange (by nearly any metric) is a mechanism known as the constant function
market maker, or CFMM. A CFMM is a particular type of DEX which allows
anyone to propose a trade (\eg, trading some amount of one asset for
another). The trade is accepted if a simple rule, which we describe
later in~\S\ref{sec:cfmms}, is met.

The prevalence of CFMMs on blockchains naturally leads to questions about
routing trades across networks or aggregations of CFMMs. For instance, suppose
that one wants to trade some amount of asset A for the greatest possible amount
of asset B. There could be many `routes' that provide this trade. For example,
we may trade asset A for asset C, and only then trade asset C for asset B.
This routing problem can be formulated as an optimization problem over the set
of CFMMs available to the user for trading. Angeris et al.~\cite{angeris2022optimal}
showed that the general problem of routing is a convex program for concave
utilities, ignoring blockchain transactions costs, though special cases of the
routing problem have been studied previously~\cite{wang2022, danos2021}.

\paragraph{This paper.} In this paper, we apply a decomposition method to the
optimal routing problem, which results in an algorithm that easily parallelizes
across all DEXs. To solve the subproblems of the algorithm, we formalize the
notions of swap markets, bounded liquidity, and aggregate CFMMs (such as
Uniswap v3) and discuss their properties. Finally, we demonstrate that our
algorithm for optimal routing is efficient, practical, and can handle the large
variety of CFMMs that exist on chain today.

\section{Optimal routing}\label{sec:opt_route}
In this section, we define the general problem of optimal routing and give concrete
examples along with some basic properties.

\paragraph{Assets.} In the optimal routing problem, we have a global labeling
of $n$ assets which we are allowed to trade, indexed by $j=1, \dots, n$
throughout this paper. We will sometimes refer to this `global collection' as
the \emph{universe} of assets that we can trade.

\paragraph{Trading sets.} Additionally, in this problem, we have a number of
markets $i=1, \dots, m$ (usually constant function market makers, or
collections thereof, which we discuss in~\S\ref{sec:cfmms}) which trade a
subset of the universe of tokens of size $n_i$. We define market $i$'s
behavior, at the time of the trade, via its \emph{trading set} $T_i \subseteq
\reals^{n_i}$. This trading set behaves in the following way: any trader is
able to propose a \emph{trade} consisting of a basket of assets $\Delta_i \in
\reals^{n_i}$, where positive entries of $\Delta_i$ denote that the trader
receives those tokens from the market, while negative values denote that the
trader tenders those tokens to the market. 
(Note that the baskets here are of a subset of
the universe of tokens which the market trades.) The market then \emph{accepts}
this trade (\ie, takes the negative elements in $\Delta_i$ from the trader and 
gives the positive elements in $\Delta_i$ to the trader) whenever
\[
    \Delta_i \in T_i.
\]
We make two assumptions about the sets $T_i$. One, that the set $T_i$ is a
closed convex set, and, two, that the zero trade is always an acceptable trade,
\ie, $0 \in T_i$. All existing DEXs that are known to the authors have a trading set that satisfies these
conditions.

\paragraph{Local and global indexing.} Each market $i$ trades only a subset
of $n_i$ tokens from the universe of tokens, so we introduce the matrices
$A_i \in \reals^{n \times n_i}$ to connect the local indices to
the global indices. These matrices are defined such that
$A_i\Delta_i$ yields the total amount of assets the trader tendered or
received from market $i$, in the global indices. For example,
if our universe has 3 tokens and market $i$ trades the tokens 2 and 3,
then
\[
    A_i = \bmat{0 & 0 \\ 1 & 0 \\ 0 & 1 }.
\]
Written another way, $(A_i)_{jk} = 1$ if token $k$ in the market's local index
corresponds to global token index $j$, and $(A_i)_{jk}=0$ otherwise. We note
that the ordering of tokens in the local index does not need to be the same as
the global ordering. 

\paragraph{Network trade vector.}
By summing the net trade in each market, after mapping the local indices to the
global indices, we obtain the \emph{network trade vector}
\[
    \Psi = \sum_{i=1}^m A_i\Delta_i.
\]
We can interpret $\Psi$ as the net trade across the network of all markets. 
If $\Psi_i > 0$, we receive some amount of asset $i$ after executing all trades 
$\left\{\Delta_i\right\}_{i=1}^m$. On the other hand, if $\Psi_i < 0$,
we tender some of asset $i$ to the network. Note that having $\Psi_i = 0$ does
not imply we do not trade asset $i$; it only means that, after executing all
trades, we received as much as we tendered.

\paragraph{Network trade utility.}
Now that we have defined the network trade vector, we introduce a utility
function $U:\reals^n \to \reals \cup \{-\infty\}$ that gives the trader's utility
of a net trade $\Psi$. We assume that $U$ is concave and increasing (\ie, we assume
all assets have value with potentially diminishing returns). Furthermore,
we will use infinite values of $U$ to encode constraints; a trade $\Psi$ such
that $U(\Psi) = -\infty$ is unacceptable to the trader. We can choose $U$ to
encode several important actions in markets, including liquidating or purchasing
a basket of assets and finding arbitrage. See~\cite[\S 5.2]{angerisConstantFunctionMarket2021}
for several examples.

\paragraph{Optimal routing problem.}
The \emph{optimal routing problem} is then the problem of finding a set of
valid trades that maximizes the trader's utility:
\begin{equation}\label{eq:opt-routing-general}
\begin{aligned}
& \text{maximize}      && U(\Psi) \\
& \text{subject to}    && \Psi = \sum_{i=1}^m A_i\Delta_i \\
                      &&& \Delta_i \in T_i, \qquad i = 1, \dots, m.
\end{aligned}
\end{equation}
The problem variables are the network trade vector $\Psi \in \reals^n$ and
trades with each market $\Delta_i \in \reals^{n_i}$, while problem data are the
utility function $U: \reals^n \to \reals \cup \{\infty\}$, the matrices $A_i
\in \reals^{n \times n_i}$, and the trading sets $T_i\subseteq \reals^{n_i}$,
where $i=1, \dots, m$. Since the trading sets are convex and the utility
function is concave, this problem is a convex optimization problem. In the
subsequent sections, we will use basic results of convex optimization to
construct an efficient algorithm to solve
problem~\eqref{eq:opt-routing-general}.

\subsection{Constant function market makers}\label{sec:cfmms}
Most decentralized exchanges, such as Uniswap~v2, Balancer, Curve, among
others, are currently organized as \emph{constant function market makers}
(CFMMs, for short) or collections of CFMMs (such as
Uniswap~v3)~\cite{angerisImprovedPriceOracles2020,angerisConstantFunctionMarket2021}.
A constant function market maker is a type of permissionless market that allows
anyone to trade baskets of, say, $r$, assets for other baskets of these same
$s$ assets, subject to a simple set of rules which we describe below.

\paragraph{Reserves and trading functions.} A constant function market maker,
which allows $r$ tokens to be traded, is defined by two properties: its
\emph{reserves} $R \in \reals_+^r$, where $R_j$ denotes the amount of asset $j$
available to the CFMM, and a \emph{trading function} which is a concave function $\phi:
\reals_+^r \to \reals$, which specifies the CFMM's behavior and its
\emph{trading fee} $0 < \gamma \le 1$.

\paragraph{Acceptance condition.} Any user is allowed to submit a trade to a
CFMM, which is, from before, a vector $\Delta \in \reals^r$. The submitted
trade is then accepted if the following condition holds:
\begin{equation}\label{eq:acceptance-condition}
    \phi(R - \gamma\Delta_- - \Delta_+) \ge \phi(R),
\end{equation}
and $R - \gamma\Delta_- - \Delta_+ \ge 0$. Here, we denote $\Delta_+$ to be the
`elementwise positive part' of $\Delta$, \ie, $(\Delta_+)_j = \max\{\Delta_j,
0\}$ and $\Delta_-$ to be the `elementwise negative part' of $\Delta$, \ie,
$(\Delta_-)_j = \min\{\Delta_j, 0\}$ for every asset $j=1, \dots, r$. The
basket of assets $\Delta_+$ may sometimes be called the `received basket' and
$\Delta_-$ may sometimes be called the `tendered basket' (see, \eg,~\cite{angerisConstantFunctionMarket2021}).
Note that the trading set $T$, for a CFMM, is exactly the set of $\Delta$ such
that~\eqref{eq:acceptance-condition} holds,
\begin{equation}\label{eq:trading-set}
    T = \{\Delta \in \reals^r \mid \phi(R - \gamma\Delta_- - \Delta_+) \ge \phi(R)\}.
\end{equation}
It is clear that $0 \in T$, and it is not difficult to show that $T$ is convex
whenever $\phi$ is concave, which is true for all trading functions used in
practice. If the trade is accepted then the CFMM pays out $\Delta_+$ from its
reserves and receives $-\Delta_-$ from the trader, which means the reserves are
updated in the following way:
\[
    R \gets R - \Delta_- - \Delta_+.
\]
The acceptance condition~\eqref{eq:acceptance-condition} can then be
interpreted as: the CFMM accepts a trade only when its trading function,
evaluated on the `post-trade' reserves with the tendered basket discounted by
$\gamma$, is at least as large as its value when evaluated on the current
reserves. 

It can be additionally shown that the trade acceptance
conditions in terms of the trading function $\phi$ and in terms of the trading
set $T$ are equivalent in the sense that every trading set has a function
$\phi$ which generates it~\cite{angerisImprovedPriceOracles2020}, under some basic
conditions.

\paragraph{Examples.}
Almost all examples of decentralized exchanges currently in production are
constant function market makers. For example, the most popular trading function
(as measured by most metrics) is the product trading function:
\[
    \phi(R) = \sqrt{R_1R_2},
\]
originally proposed for Uniswap~\cite{uniswap} and a `bounded liquidity' variation of this function:
\begin{equation}\label{eq:bounded-liquidity}
    \phi(R) = \sqrt{(R_1 + \alpha)(R_2 + \beta)},
\end{equation}
used in Uniswap v3~\cite{uniswapv3}, with $\alpha, \beta \ge 0$. 
Other examples include the weighted geometric mean (as used by Balancer~\cite{balancer})
\begin{equation}\label{eq:geom}
    \phi(R) = \prod_{i=1}^r R_i^{w_i},
\end{equation}
where $r$ is the number of assets the exchange trades, and $w \in \reals_+^r$
with $\ones^Tw = 1$ are known as the weights, along with the Curve trading
function
\[
    \phi(R) = \alpha \ones^TR - \left(\prod_{i=1}^r R_i^{-1}\right),
\]
where $\alpha > 0$ is a parameter set by the CFMM~\cite{egorovStableSwapEfficientMechanism}.
Note that the `product' trading
function is the special case of the weighted geometric mean function when $r=2$ and
$w_1 = w_2 = 1/2$.

\paragraph{Aggregate CFMMs.} \label{sec:collection}
In some special cases, such as in Uniswap v3, it
is reasonable to consider an \emph{aggregate CFMM}, which we define as a
{collection} of CFMMs, which all trade the same assets, as part of a
single `big' trading set. A specific instance of an aggregate CFMM currently
used in practice is in Uniswap v3~\cite{uniswapv3}. Any `pool' in this exchange
is actually a collection of CFMMs with the `bounded liquidity' variation of the
product trading function, shown in~\eqref{eq:bounded-liquidity}. We will see
that we can treat these `aggregate CFMMs' in a special way in order to significantly
improve performance.

\section{An efficient algorithm}\label{sec:eff_alg}
A common way of solving problems such as
problem~\eqref{eq:opt-routing-general}, where we have a set of variables
coupled by only a single constraint, is to use a decomposition
method~\cite{dantzig1960decomposition, bertsekasNonlinearProgramming2016}. 
The general idea of these methods is to solve the original problem by splitting
it into a sequence of easy subproblems that can be solved independently.
In this section, we will see that applying a decomposition
method to the optimal routing problem gives a solution method which parallelizes
over all markets. Furthermore, it gives a clean programmatic interface; we only need
to be able to find arbitrage for a market, given a set of reference prices. This interface allows us
to more easily include a number of important decentralized exchanges, such as
Uniswap v3.

\subsection{Dual decomposition}
To apply the dual decomposition method, we
first take the coupling constraint of problem~\eqref{eq:opt-routing-general},
\[
\Psi = \sum_{i=1}^m A_i\Delta_i,
\]
and relax it to a linear penalty in the objective, parametrized by some vector
$\nu \in \reals^n$. (We will show in~\S \ref{sec:dual-problem} that the only reasonable choice of
$\nu$ is a market clearing price, sometimes called a no-arbitrage price, and that
this choice actually results in a relaxation that is tight; \ie, a solution for
this relaxation also satisfies the original coupling constraint.) This relaxation results
in the following problem:
\[
\begin{aligned}
    & \text{maximize}      && {\textstyle U(\Psi) - \nu^T(\Psi - \sum_{i=1}^m A_i\Delta_i)}\\
    & \text{subject to}    && \Delta_i \in T_i, \quad i=1, \dots, m,
\end{aligned}
\]
where the variables are the network trade vector $\Psi \in \reals^n$ and the
trades are $\Delta_i \in \reals^{n_i}$ for each market $i=1,
\dots, m$. Note that this formulation can be viewed as a family of problems parametrized
by the vector $\nu$.

A simple observation is that this new problem is actually separable over
all of its variables. We can see this by rearranging the objective:
\begin{equation}\label{eq:rewriting}
\begin{aligned}
    & \text{maximize}      && U(\Psi) - \nu^T\Psi + {\textstyle \sum_{i=1}^m (A_i^T\nu)^T\Delta_i}\\
    & \text{subject to}    && \Delta_i \in T_i, \quad i=1, \dots, m.
\end{aligned}
\end{equation}
Since there are no additional coupling constraints, we can solve for $\Psi$
and each of the $\Delta_i$ with $i=1, \dots, m$ separately.

\paragraph{Subproblems.} This method gives two types of subproblems, each
depending on $\nu$. The first, over $\Psi$, is relatively simple:
\begin{equation}\label{eq:first}
    \begin{aligned}
        & \text{maximize} && U(\Psi) - \nu^T\Psi,
    \end{aligned}
\end{equation}
and can be recognized as a slightly transformed version of the Fenchel conjugate~\cite[\S 3.3]{cvxbook}.
We will write its optimal value (which depends on $\nu$) as
\[
    \bar U(\nu) = \sup_{\Psi} \left(U(\Psi) - \nu^T\Psi\right).
\]
The function $\bar U$ can be easily derived in closed form for a number of
functions $U$. Additionally, since $\bar U$ is a supremum over an affine family
of functions parametrized by $\nu$, it is a convex function of
$\nu$~\cite[\S 3.2.3]{cvxbook}. (We will use this fact soon.) Another important thing to note
is that unless $\nu \ge 0$, the function $\bar U(\nu)$ will evaluate to $+\infty$. 
This can be interpreted as an implicit constraint on $\nu$.

The second type of problem is over each trade $\Delta_i$ for $i=1, \dots, m$,
and can be written, for each market $i$, as
\begin{equation}\label{eq:arb}
    \begin{aligned}
        & \text{maximize} && (A_i^T\nu)^T\Delta_i\\
        & \text{subject to} && \Delta_i \in T_i.
    \end{aligned}
\end{equation}
We will write its optimal value, which depends on $A_i^T\nu$, as
$\arb_i(A_i^T\nu)$. Problem~\eqref{eq:arb} can be recognized as the \emph{optimal
arbitrage problem} (see, \eg,~\cite{angerisConstantFunctionMarket2021}) for market $i$, when the external
market price, or reference market price, is equal
to $A_i^T\nu$. Since $\arb_i(A_i^T\nu)$ is also defined as a supremum over a family of affine
functions of $\nu$, it too is a convex function of $\nu$. Solutions to the
optimal arbitrage problem are known, in closed form, for a number of trading functions.
(See appendix~\ref{app:closed-form} for some examples.)

\paragraph{Dual variables as prices.} The optimal solution to
problem~\eqref{eq:arb}, given by $\Delta_i^\star$, is a point
$\Delta_i^\star$ in $T_i$ such that there exists a supporting hyperplane to the set $T_i$ at
$\Delta_i^\star$ with slope $A_i^T\nu$~\cite[\S 5.6]{cvxbook}. We can interpret
these slopes as the `marginal prices' of the $n_i$ assets, since, letting
$\delta \in \reals^{n_i}$ be a small deviation from the trade $\Delta_i^\star$,
we have, writing $\tilde \nu = A_i^T\nu$ as the weights of $\nu$ in the local
indexing:
\[
    \tilde \nu^T(\Delta_i^\star + \delta) \le \tilde\nu^T\Delta_i^\star,
\]
for every $\delta$ with $\Delta_i^\star + \delta \in T_i$. (By definition of optimality.)
Canceling terms, we find:
\[
    \tilde \nu^T\delta \le 0.
\]
If, for example, $\delta_i$ and $\delta_j$ are the only two nonzero entries of
$\delta$, we would have
\[
    \delta_i \le -\frac{\tilde\nu_j}{\tilde \nu_i}\delta_j,
\]
so the exchange rate between $i$ and $j$ is at most $\tilde\nu_i/\tilde\nu_j$.
This observation lets us interpret the dual variables $\tilde\nu$ (and
therefore the dual variables $\nu$) as `marginal prices', up to a constant multiple.

\subsection{The dual problem}\label{sec:dual-problem}
The objective value of problem~\eqref{eq:rewriting}, which is a function of
$\nu$, can then be written as
\begin{equation}\label{eq:dual-function}
    g(\nu) = \bar U(\nu) + \sum_{i=1}^m \arb_i(A_i^T\nu).
\end{equation}
This function $g: \reals^n \to \reals$ is called the \emph{dual function}.
Since $g$ is the sum of convex functions, it too is convex. The \emph{dual problem}
is the problem of minimizing the dual function,
\begin{equation}\label{eq:dual}
    \begin{aligned}
        & \text{minimize} && g(\nu),
    \end{aligned}
\end{equation}
over the dual variable $\nu \in \reals^n$, which is a convex optimization
problem since $g$ is a convex function.

\paragraph{Dual optimality.} While we have defined the dual problem, we
have not discussed how it relates to the original routing problem we are
attempting to solve, problem~\eqref{eq:opt-routing-general}. Let $\nu^\star$ be
a solution to the dual problem~\eqref{eq:dual}. Assuming that the dual function
is differentiable at $\nu^\star$, the first order, unconstrained optimality
conditions for problem~\eqref{eq:dual} are that
\[
    \nabla g(\nu^\star) = 0.
\]
(The function $g$ need not be differentiable, in which case a similar, 
but more careful, argument holds using subgradient calculus.)
It is not hard to show that if $\bar U$ is differentiable at $\nu^\star$, then
its gradient must be $\nabla \bar U(\nu^\star) = -\Psi^\star$, where
$\Psi^\star$ is the solution to the first subproblem~\eqref{eq:first}, with
$\nu^\star$. (This follows from the fact that the gradient of a maximum, when
differentiable, is the gradient of the argmax.) Similarly, the gradient of
$\arb_i$ when evaluated at $A_i^T\nu^\star$ is $ \Delta_i^\star$, where
$\Delta_i^\star$ is a solution to problem~\eqref{eq:arb}
with marginal prices $A_i^T\nu^\star$, for each market $i=1, \dots, m$. Using
the chain rule, we then have:
\begin{equation}\label{eq:grad}
    0 = \nabla g(\nu^\star) = -\Psi^\star + \sum_{i=1}^m
        A_i\Delta_i^\star.
\end{equation}
Note that this is exactly the coupling constraint of
problem~\eqref{eq:opt-routing-general}. In other words, when the linear
penalties $\nu^\star$ are chosen optimally (\ie, chosen such that they minimize
the dual problem~\eqref{eq:dual}) then the optimal solutions for
subproblems~\eqref{eq:first} and~\eqref{eq:arb} automatically satisfy the
coupling constraint. Because problem~\eqref{eq:rewriting} is a relaxation of
the original problem~\eqref{eq:opt-routing-general} for any choice of $\nu$,
any solution to problem~\eqref{eq:rewriting} that satisfies the coupling
constraint of problem~\eqref{eq:opt-routing-general} must also be a solution to
this original problem. All that remains is the question of finding a solution
$\nu^\star$ to the dual problem~\eqref{eq:dual}.

\subsection{Solving the dual problem}
The dual problem~\eqref{eq:dual} is a convex optimization problem that is easily
solvable in practice, even for very large $n$ and $m$. In many cases, we can use a
number of off-the-shelf solvers such as SCS~\cite{odonoghueConicOptimizationOperator2016},
Hypatia~\cite{coey2021solving}, 
and Mosek~\cite{mosek}. For example, a relatively
effective way of minimizing functions when the gradient is easily evaluated is the
L-BFGS-B algorithm~\cite{byrd1995limited,zhu1997algorithm,morales2011remark}: 
given a way of evaluating the dual
function $g(\nu)$ and its gradient $\nabla g(\nu)$ at some point $\nu$, the
algorithm will find an optimal $\nu^\star$ fairly quickly in practice. (See~\S
\ref{sec:performance} for timings.) By definition, the function $g$ is easy to evaluate if the
subproblems~\eqref{eq:first} and~\eqref{eq:arb} are easy to evaluate.
Additionally the right hand side of equation~\eqref{eq:grad} gives us a way of
evaluating the gradient $\nabla g$, essentially for free,
since we typically receive the optimal $\Psi^\star$ and
$\Delta_i^\star$ as a consequence of computing $\bar U$ and $\arb_i$.

\paragraph{Interface.}\label{sec:interface}
In order for a user to specify and solve the dual
problem~\eqref{eq:dual} (and therefore the original problem) it suffices for
the user to specify (a) some way of evaluating $\bar U$ and its optimal $\Psi$
for problem~\eqref{eq:first} and (b) some way of evaluating the arbitrage
problem~\eqref{eq:arb} and its optimal trade $\Delta_i^\star$ for each
market $i$ that the user wishes to include. New markets can be easily added by
simply specifying how to arbitrage them, which, as we will see next, turns out
to be straightforward for most practical decentralized exchanges. The Julia interface
required for the software package described in~\S\ref{sec:implementation}
is a concretization of the interface described here.

\section{Swap markets}\label{sec:swap}
In practice, most markets trade only two assets; we will refer to these kinds
of markets as \emph{swap markets}. Because these markets are so common, the
performance of our algorithm is primarily governed by its ability to
solve~\eqref{eq:arb} quickly on these two asset markets. We show practical
examples of these computations in appendix~\ref{app:closed-form}. In this section,
we will suppress the index $i$ with the understanding that we are referring to a specific
market $i$.

\subsection{General swap markets}
Swap markets are simple to deal with because their trading behavior is completely
specified by the \emph{forward exchange function}~\cite{angerisConstantFunctionMarket2021}
for each of the two assets.
In what follows, the forward trading function $f_1$ will denote the maximum amount of
asset 2 that can be received by trading some fixed amount $\delta_1$ of asset 1, \ie,
if $T \subseteq \reals^2$ is the trading set for a specific swap market,
then
\[
    f_1(\delta_1) = \sup \{ \lambda_2 \mid (-\delta_1, \lambda_2) \in T\}, \quad  f_2(\delta_2) = \sup \{ \lambda_1 \mid (\lambda_1, -\delta_2) \in T\}.
\]
In other words, $f_1(\delta_1)$ is defined as the largest amount $\lambda_2$ of
token 2 that one can receive for tendering a basket of $(\delta_1, 0)$ to the
market. The forward trading function $f_2$ has a similar interpretation. If $f_1(\delta_1)$
is finite, then this supremum is achieved since the set $T$ is closed.

\paragraph{Trading function.}
If the set $T$ has a simple trading function representation, as
in~\eqref{eq:trading-set}, it is not hard to show that the function $f_1$ is
the unique (pointwise largest) function that satisfies
\begin{equation}\label{eq:func-form}
\phi(R_1 + \gamma\delta_1, R_2 - f_1(\delta_1)) = \phi(R_1, R_2).
\end{equation}
whenever $\phi$ is nondecreasing, which may be assumed for all CFMMs~\cite{angerisImprovedPriceOracles2020},
and similarly for $f_2$. (Note the equality here, compared to the inequality in
the original definition~\eqref{eq:acceptance-condition}.) 

\paragraph{Properties.}
The functions $f_1$ and $f_2$ are concave, since the trading set
$T$ is convex, and nonnegative, since $0 \in T$ by assumption. Additionally, we
can interpret the directional derivative of $f_j$ as the current marginal price
of the received asset, denominated in the tendered asset. Specifically, we
define
\begin{equation}\label{eq:marginal}
    f'_j(\delta_j) = \lim_{h \to 0^+} \frac{f_j(\delta_j + h) - f_j(\delta_j)}{h}.
\end{equation}
This derivative is sometimes referred to as the price impact
function~\cite{angeris2022does}. Intuitively, $f_1'(0)$ is the current price of
asset $1$ quoted by the swap market before any trade is made, and
$f_1'(\delta)$ is the price quoted by the market to add an additional $\eps$
units of asset $1$ to a trade of size $\delta$, for very small $\eps$. We note
that in the presence of fees, the marginal price to add to a trade of size
$\delta$, \ie, $f_1'(\delta)$, will be lower than the price to do so after the
trade has been made~\cite{angerisImprovedPriceOracles2020}.

\paragraph{Swap market arbitrage problem.}
Equipped with the forward exchange function, we can specialize~\eqref{eq:arb}.
Overloading notation slightly by writing $(\nu_1, \nu_2) \ge 0$ for $A_i^T\nu$ we
define the swap market arbitrage problem for a market with forward exchange
function $f_1$:
\begin{equation}\label{eq:scalar-arb}
\begin{aligned}
& \text{maximize}   && -\nu_1\delta_1 + \nu_2f_1(\delta_1) \\
& \text{subject to} && \delta_1 \ge 0,
\end{aligned}
\end{equation}
with variable $\delta_1 \in \reals$
We can also define a similar arbitrage problem for $f_2$:
\[
\begin{aligned}
& \text{maximize}   && \nu_1 f_2(\delta_2) - \nu_2\delta_2\\
& \text{subject to} && \delta_2 \ge 0,
\end{aligned}
\]
with variable $\delta_2 \in \reals$. Since $f_1$ and $f_2$ are concave, both
problems are evidently convex optimization problems of one variable. Because
they are scalar problems, these problems can be easily solved by bisection or
ternary search. The final solution is to take whichever of these two problems
has the largest objective value and return the pair in the correct order. For
example, if the first problem~\eqref{eq:scalar-arb} has the highest objective
value with a solution $\delta_1^\star$, then $\Delta^\star = (-\delta_1^\star,
f(\delta_1^\star))$ is a solution to the original arbitrage
problem~\eqref{eq:arb}. (For many practical trading sets $T$, it can be shown that at
most one problem will have strictly positive objective value, so it is possible
to `short-circuit' solving both problems if the first evaluation has positive optimal value.)

\paragraph{Problem properties.} One way to view each of these problems is that
they `separate' the solution space of the original arbitrage
problem~\eqref{eq:arb} into two cases: one where an optimal solution
$\Delta^\star$ for~\eqref{eq:arb} has $\Delta_1^\star \le 0$ and one where an
optimal solution has $\Delta_2^\star \le 0$. (Any optimal point $\Delta^\star$
for the original arbitrage problem~\eqref{eq:arb} will never have both
$\Delta^\star_1 < 0$ and $\Delta^\star_2 < 0$ as that would be strictly worse
than the 0 trade for $\nu > 0$, and no reasonable market will have
$\Delta^\star_1 > 0$ and $\Delta^\star_2 > 0$ since the market would be
otherwise `tendering free money' to the trader.) This observation means that, in order to
find an optimal solution to the original optimal arbitrage
problem~\eqref{eq:arb}, it suffices to solve two scalar convex optimization
problems.

\paragraph{Optimality conditions.}
The optimality conditions for problem~\eqref{eq:scalar-arb} are
that, if
\begin{equation}\label{eq:swap-opt-cond}
    \nu_2 f_1'(0) \le \nu_1
\end{equation}
then $\delta_1^\star = 0$ is a solution. Otherwise, we have
\[
    \delta^\star_1 = \sup \{\delta \ge 0 \mid \nu_2f_1'(\delta) \ge \nu_1\}.
\]
Similar conditions hold for the problem over $\delta_2$. If the function $f_1'$
is continuous, not just semicontinuous, then the expression above simplifies to
finding a root of a monotone function:
\begin{equation}\label{eq:root}
\nu_2 f_1'(\delta^\star_1) = \nu_1.
\end{equation}
If there is no root and condition~\eqref{eq:swap-opt-cond} does not hold, then
$\delta_1^\star = \infty$. However, the solution will be finite for any trading
set that does not contain a line, \ie, the market does not have
`infinite liquidity' at a specific price.

\paragraph{No-trade condition.}
Note that using the inequality~\eqref{eq:swap-opt-cond} gives us a simple way
of verifying whether we will make any trade with market $T$, given some prices
$\nu_1$ and $\nu_2$. In particular, the zero trade is optimal whenever
\[
    f_1'(0) \le \frac{\nu_1}{\nu_2} \le \frac{1}{f_2'(0)}.
\]
We can view the interval $[f_1'(0), 1/f_2'(0)]$ as a type of `bid-ask spread'
for the market with trading set $T$. (In constant function market makers, this
spread corresponds to the fee $\gamma$ taken from the trader.) This `no-trade condition' lets us
save potentially wasted effort of computing an optimal arbitrage trade as,
in practice, most trades in the original problem will be 0.

\paragraph{Bounded liquidity.} \label{sec:bounded-liq}
In some cases, we can
easily check not only when a trade will not be made (say, using condition~\eqref{eq:swap-opt-cond}),
but also  when the `largest possible trade' will be made. (We will define
what this means next.) Markets for which there is a `largest possible trade'
are called bounded liquidity markets. We say a market has \emph{bounded
liquidity in asset 2} if there is a finite $\delta_1$ such that $f_1(\delta_1)
= \sup f_1$, and similarly for $f_2$. In other words, there is a finite input
$\delta_1$ which will give the maximum possible amount of asset 2 out. A market
has \emph{bounded liquidity} if it has bounded liquidity on both of its assets.
A bounded liquidity market then has a notion of a `minimum price'. First,
define
\[
    \delta_1^- = \inf \{\delta_1 \ge 0 \mid f_1(\delta_1) = \sup f_1\},
\]
\ie, $\delta_1^-$ is the smallest amount of asset 1 that can be tendered to receive
the maximum amount the market is able to supply. We can then define the \emph{minimum
supported price} as the left derivative of $f_1$ at $\delta_1^-$:
\[
    f_1^-(\delta_1^-) = \lim_{h \to 0^+} \frac{f(\delta_1^-) - f(\delta_1^- - h)}{h}.
\]
The first-order optimality conditions imply that $\delta_1^-$ is a solution to the scalar optimal arbitrage
problem~\eqref{eq:scalar-arb} whenever
\[
    f_1^-(\delta_1^-) \ge \frac{\nu_1}{\nu_2}.
\]
In English, this can be stated as: if the minimum supported marginal price we
receive for $\delta_1^-$ is still larger than the price being arbitraged against,
$\nu_1/\nu_2$, it is optimal to take all available liquidity from the
market. Using the same definitions for $f_2$, we find that the only time
the full problem~\eqref{eq:scalar-arb} needs to be solved is when
the price being arbitraged against $\nu_1/\nu_2$ lies in the interval
\begin{equation}\label{eq:interval}
    f_1^-(\delta_1^-) < \frac{\nu_1}{\nu_2} < \frac{1}{f_2^-(\delta_2^-)}.
\end{equation}
(It may be the case that $f_2^-(\delta_2^-) = 0$ in which case we define
the right hand side to be $\infty$.) We will call this interval of prices
the \emph{active interval} for a bounded liquidity market.

\paragraph{Example.}\label{sec:bounded}
In the case of Uniswap v3~\cite{uniswapv3}, we have a collection of, say, $i=1,
\dots, s$ bounded liquidity product functions~\eqref{eq:bounded-liquidity},
where the parameters $\alpha_k, \beta_k > 0$ are chosen such that all of the
active price intervals, as defined in~\eqref{eq:interval}, are disjoint. (An
explicit form for this trading function is given in the appendix,
equation~\eqref{eq:e-interval}.) Solving the arbitrage
problem~\eqref{eq:scalar-arb} over this collection of CFMMs is relatively
simple. Since all of the intervals are disjoint, any price $\nu_1/\nu_2$ can
lie in at most one of the active intervals. We therefore do not need to compute the
optimal trade for any interval, except the single interval where $\nu_1/\nu_2$
lies, which can be done in closed form.
We also note that this `trick' applies to any
collection of bounded liquidity markets with disjoint active price intervals.

\section{Implementation}
\label{sec:implementation}
We have implemented this algorithm in \texttt{CFMMRouter.jl}, a Julia~\cite{bezansonJuliaFreshApproach2017} 
package for solving the optimal routing problem.
Our implementation is available at
\begin{center}
   \texttt{https://github.com/bcc-research/CFMMRouter.jl}
\end{center}
and includes implementations for both weighted geometric mean CFMMs and Uniswap 
v3. In this section, we provide a concrete Julia interface for our solver.

\subsection{Markets}
\label{sec:implementation-markets}
\paragraph{Market interface.}
As discussed in~\S\ref{sec:interface}, the only function that the
user needs to implement to solve the routing problem for a given
market is 
\begin{center}
    \jlinl{find_arb!(Δ, Λ, mkt, v)}.
\end{center}
This function solves the optimal arbitrage problem~\eqref{eq:arb}
for a market \jlinl{mkt} (which holds the relevant data about the
trading set $T$) with dual variables \jlinl{v} (corresponding
to $A_i^T\nu$ in the original problem~\eqref{eq:arb}).
It then fills the vectors \jlinl{Δ} and \jlinl{Λ} with the negative part
of the solution, $-\Delta_-^\star$, and positive part of the solution,
$\Delta_+^\star$, respectively.

For certain common markets (\eg, geometric mean and Uniswap v3), we provide
specialized, efficient implementations of \jlinl{find_arb!}. For general CFMMs
where the trading function, its gradient, and the Hessian are easy to evaluate,
one can use a general-purpose primal-dual interior point solver. For other
more complicated markets, a custom implementation may be required.

\paragraph{Swap markets.}\label{sec:swap-interface}
The discussion in~\S\ref{sec:swap} and the expression in~\eqref{eq:root}
suggests a natural, minimal interface for swap markets.
Specifically, we can define a swap market by implementing the function
\jlinl{get_price(Δ)}. This function takes in a vector of inputs \jlinl{Δ} $\in
\reals^2_+$, where we assume that only one of the two assets is being
tendered, \ie, \jlinl{Δ₁Δ₂ == 0}, and returns $f_1'(\Delta_1)$, if $\Delta_1 >
0$ or $f_2'(\Delta_2)$ if $\Delta_2 > 0$. With this \emph{price impact}
function implemented, one can use bisection to compute the solution
to~\eqref{eq:root}. When price impact function has a closed
form and is readily differentiable by hand, it is possible to use a much
faster Newton method to solve this problem. In the case where the function does
not have a simple closed form, we can use automatic differentiation (\eg, using
\texttt{ForwardDiff.jl}~\cite{RevelsLubinPapamarkou2016}) 
to generate the gradients for this function.

\paragraph{Aggregate CFMMs.}
In the special case of aggregate, bounded liquidity CFMMs, the price impact
function often does not have a closed form. On the other hand, whenever the
active price intervals are disjoint, we can use the trick presented
in~\S\ref{sec:bounded} to quickly arbitrage an aggregate CFMM. 
For example, a number of Uniswap v3
markets are actually composed of many thousands of bounded liquidity CFMMs.
Treating each of these as their own market, without any additional
considerations, significantly increases the size and solution complexity of the problem.

In this special case, each aggregate market `contains' $s$ trading sets, each
of which has disjoint active price intervals with all others. We will write
these intervals as $(p_{i}^-, p_{i}^+)$ for each trading set $i=1, \dots, s$,
and assume that these are in sorted order $p_{i-1}^+ \le p_{i}^- < p_{i}^+ \le
p_{i+1}^+$. Given some dual variables $\nu_1$ and $\nu_2$ for which to solve
the arbitrage problem~\eqref{eq:arb}, we can then run binary search over the
sorted intervals (taking $O(\log(s))$ time) to find which of the intervals the
price $\nu_1/\nu_2$ lies in. We can compute the optimal arbitrage for this
`active' trading set, and note that the remaining trading sets all have a known
optimal trade (from the discussion in~\S\ref{sec:bounded-liq}) and require only
constant time. For Uniswap v3 and other aggregate CFMMs, this algorithm is much
more efficient from both a computational and memory perspective when compared
with a direct approach that considers all $s$ trading sets separately.

\paragraph{Other functions.}
If one is solving the arbitrage problem multiple times in a row, it may be
helpful to implement the following additional functions:
\begin{enumerate}
    \item \jlinl{swap!(cfmm, Δ)}: updates \jlinl{cfmm}'s state following a trade \jlinl{Δ}.
    \item \jlinl{update_liquidity!(cfmm, [range,] L)}: adds some amount of 
    liquidity \jlinl{L} $\in \reals^2_+$, optionally includes some interval \jlinl{range = (p1, p2)}.
\end{enumerate}

\subsection{Utility functions.}\label{sec:utility-functions}
Recall that the dual problem relies on a slightly transformed version of the
Fenchel conjugate, which is the optimal value of problem~\eqref{eq:first}. To
use LBFGS-B (and most other optimization methods), we need to be able to evaluate
this function $\bar U(\nu)$ and its gradient $\nabla \bar U(\nu)$, which is the 
solution $\Psi^\star$ to~\eqref{eq:first} with parameter $\nu$.
Thus, utility functions are implemented as objects that implement the following 
interface:
\begin{itemize}
    \item \jlinl{f(objective, v)} evaluates $\bar U$ at \jlinl{v}.
    \item \jlinl{grad!(g, objective, v)} evaluates $\nabla \bar U$ at \jlinl{v} 
    and stores it in \jlinl{g}.
    \item \jlinl{lower_limit(objective)} returns the lower bound of the objective.
    \item \jlinl{upper_limit(objective)} returns the upper bound of the objective.
\end{itemize}
The lower and upper bounds can be found by deriving the conjugate function. For
example, for the `total arbitrage' objective $U(\Psi) = c^T\Psi - I(\Psi \ge
0)$, where a trader wants to tender no tokens to the network, but receive any
positive amounts out with value proportional to some nonnegative vector $c \in
\reals^n_+$, has $\bar U(\nu) = 0$ if $\nu \ge c$ and $\infty$ otherwise. Thus,
we have the bounds $c \le \nu < \infty$, and gradient $\nabla\bar U(\nu) = 0$.
We provide implementations for arbitrage and for basket
liquidations in our Julia package. (See~\cite[\S3]{angeris2022optimal} for definitions.)

\section{Numerical results}
We compare the performance of our solver against the commercial, off-the-shelf
convex optimization solver Mosek, accessed through JuMP~\cite{dunningJuMPModelingLanguage2017,legatMathOptInterfaceDataStructure2021}.
In addition, we use our solver with real, on-chain data to illustrate the
benefit of routing an order through multiple markets rather than trading with 
a single market.
Our code is available at
\begin{center}
    \texttt{https://github.com/bcc-research/router-experiments}.
\end{center}

\paragraph{Performance.}
\label{sec:performance}
We first compare the performance of our solver against Mosek~\cite{mosek}, a
widely-used, performant commercial convex optimization solver. We generate $m$
swap markets over a global universe of $2\sqrt{m}$ assets. Each market is randomly
generated with reserves uniformly sampled from the interval between $1000$ and $2000$,
denoted $R_i \sim \mathcal{U}(1000, 2000)$, and is a constant
product market with probability $0.5$ and a weighted geometric mean market with
weights $(0.8, 0.2)$ otherwise.
(These types of swap markets are common in protocols such as Balancer~\cite{balancer}.)
We run arbitrage over the set of markets, with `true prices' for each asset
randomly generated as $p_i \sim \mathcal{U}(0,1)$. For each $m$, we use the same
parameters (markets and price) for both our solver and Mosek.
Mosek is configured with default parameters.
All experiments are run on a MacBook Pro with a 2.3GHz 8-Core Intel i9 processor.
\begin{figure}
    \centering
    \includegraphics[width=0.48\linewidth]{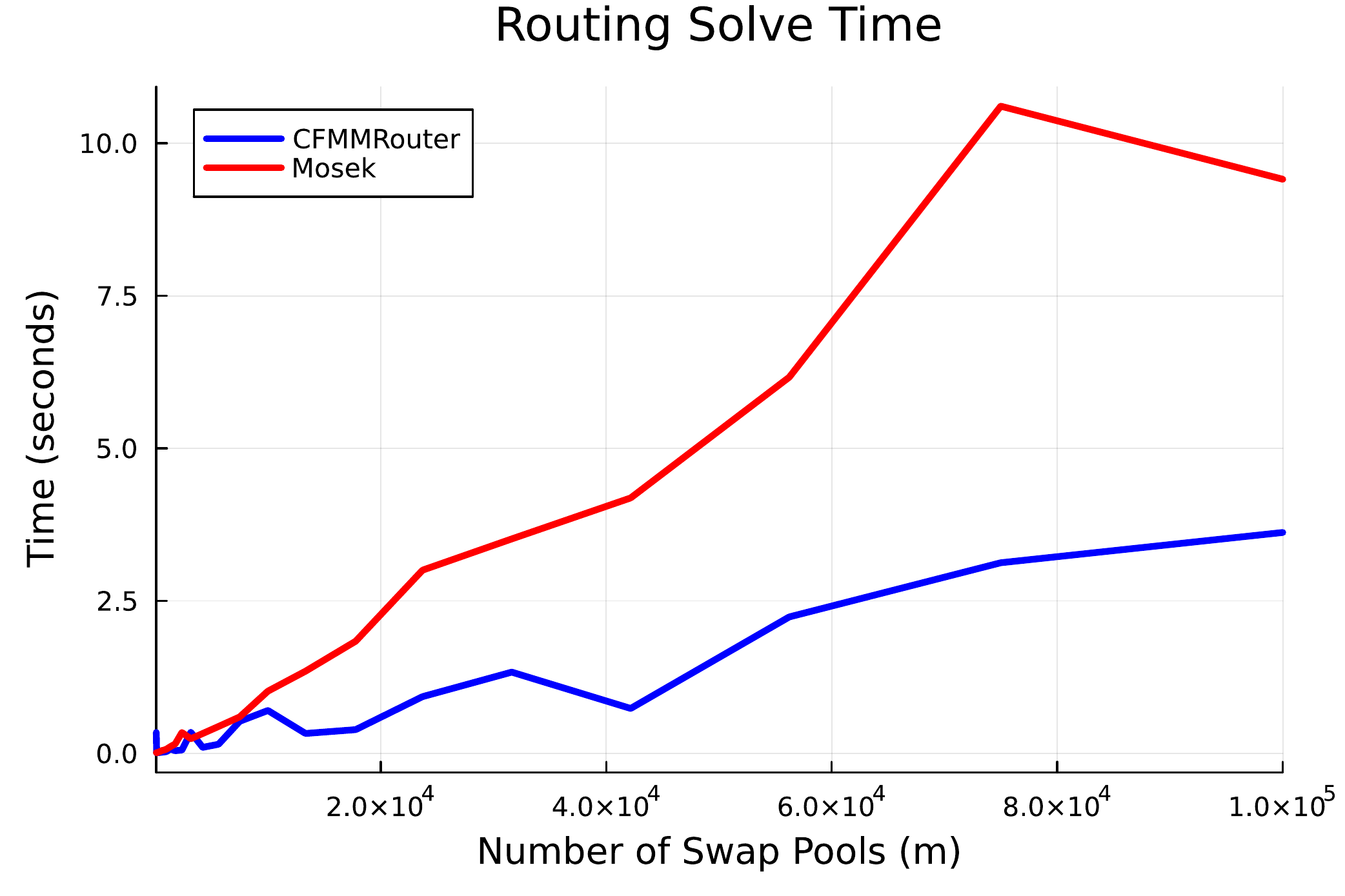}
    \includegraphics[width=0.48\linewidth]{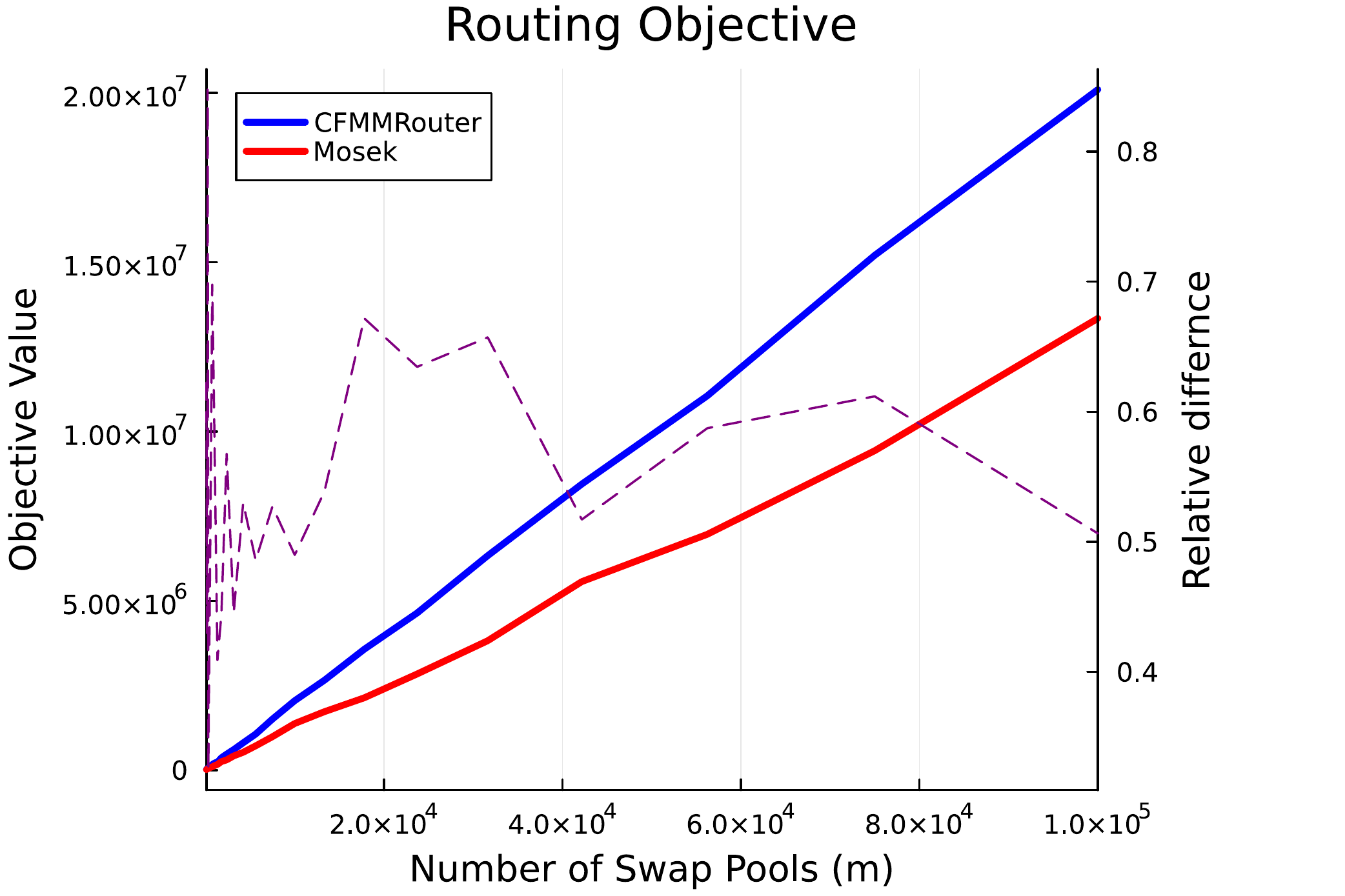}
    \caption{Solve time of Mosek vs. \texttt{CFMMRouter.jl} (left) and the resulting
    objective values for the arbitrage problem, with the dashed line indicating
    the relative increase in objective provided by our method (right).}
    \label{fig:mosek}
\end{figure}
In figure~\ref{fig:mosek}, we see that as the number of pools (and tokens) grow,
our method begins to dramatically outperform Mosek and scales quite a bit better.
We note that the weighted geometric mean markets are especially hard for Mosek,
as they must be solved as power cone constraints. Constant product markets may be
represented as second order cone constraints, which are quite a bit more efficient
for many solvers. Furthermore, our method gives a higher objective value, often by
over 50\%. We believe this increase stems from Mosek's use of an interior point 
method and numerical tolerances.
The solution returned by Mosek for each market will be strictly inside the 
associated trading set, but we know that any rational trader will choose a trade
on the boundary.

\begin{figure}
    \centering
    \includegraphics[width=0.48\linewidth]{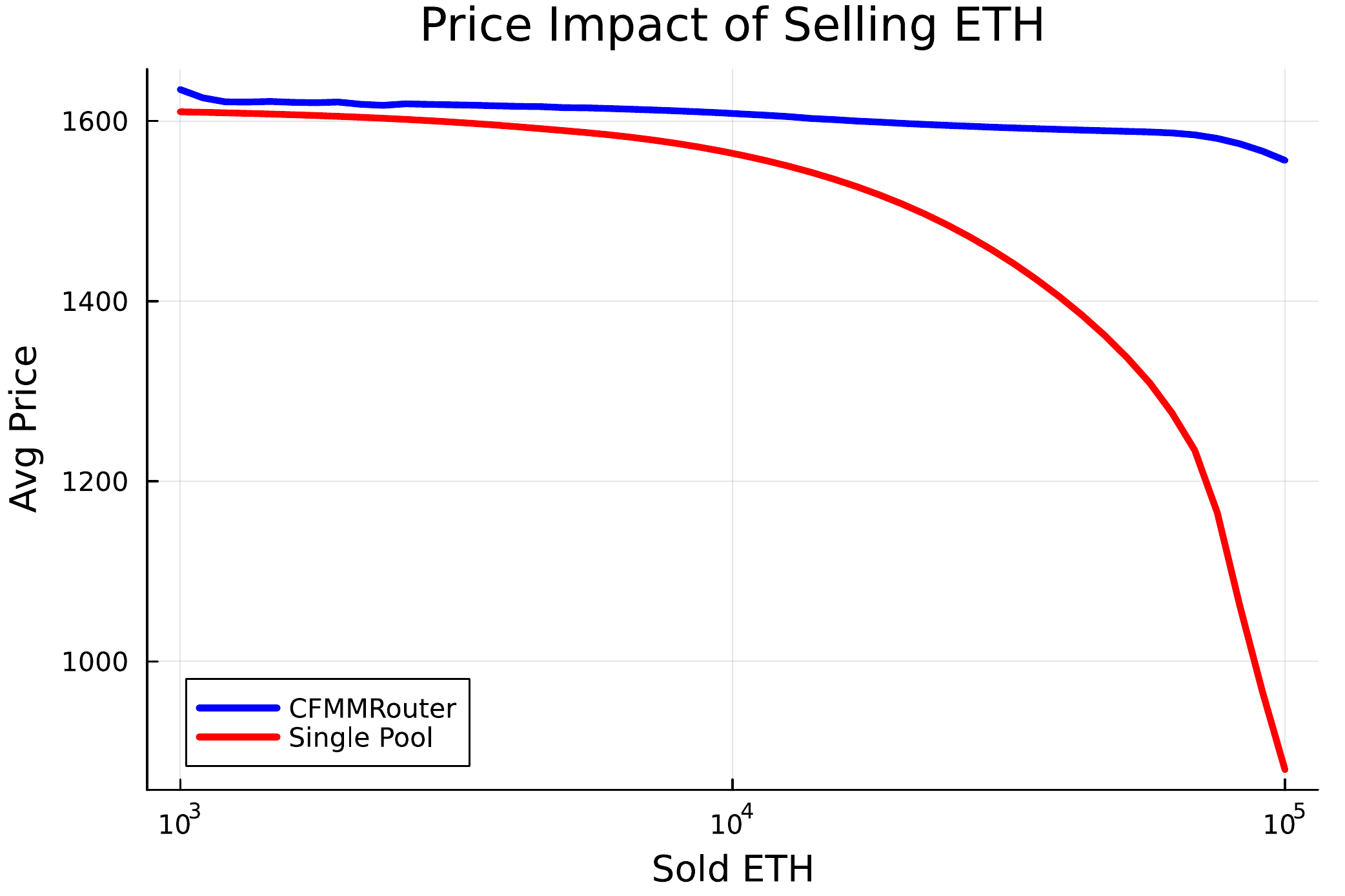}
    \includegraphics[width=0.48\linewidth]{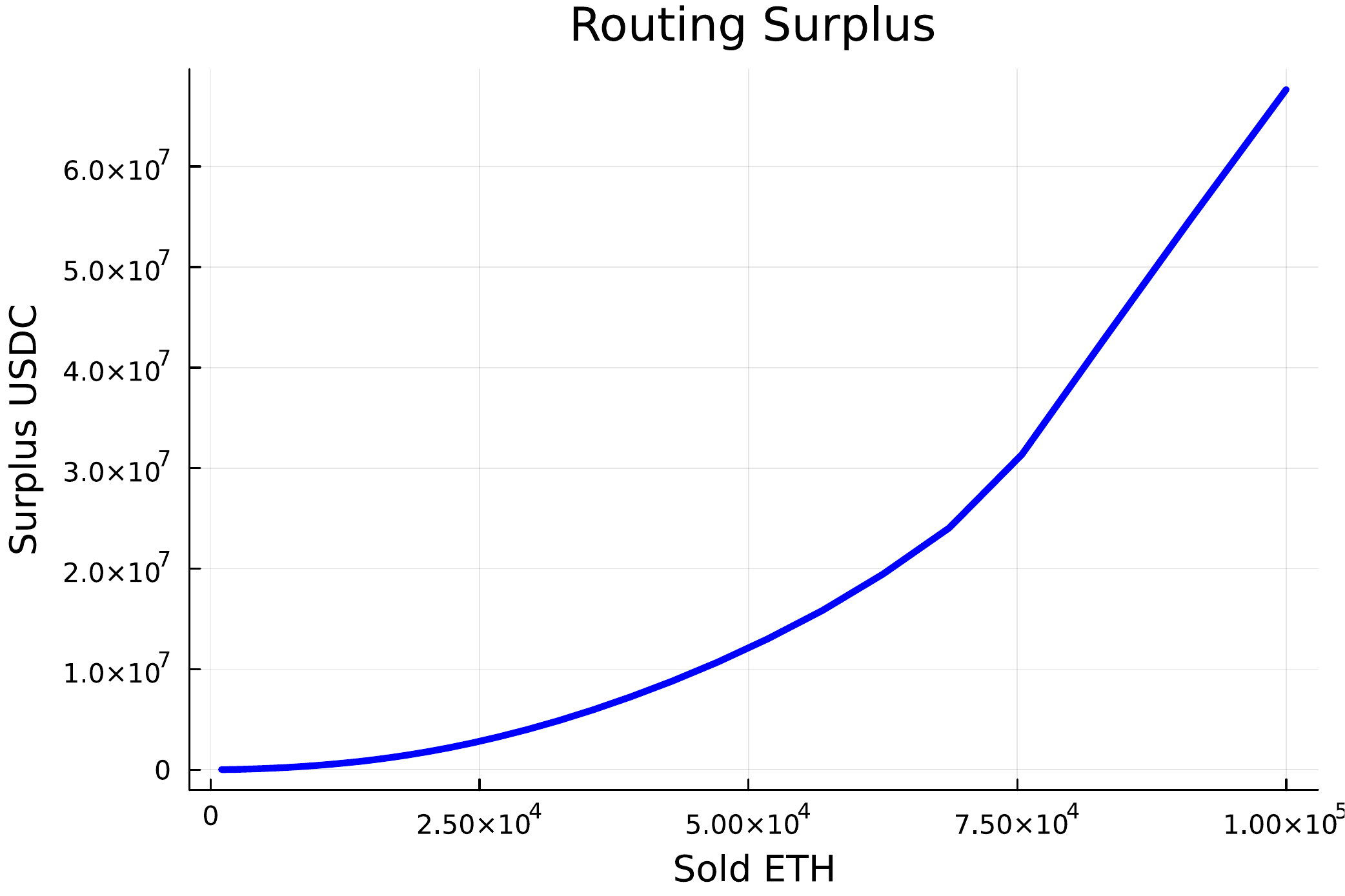}
    \caption{
    Average price of market sold ETH in routed vs.\ single-pool (left) and routed vs.\ single-pool surplus liquidation value 
    (right).}
    \label{fig:univ3}
\end{figure}

\paragraph{Real data: trading on chain.}
We show the efficacy of routing by considering a swap from WETH to USDC (\ie,
using the basket liquidation objective to sell WETH for USDC). Using on-chain
data from the end of a recent block, we show in figure~\ref{fig:univ3} 
that as the trade size increases,
routing through multiple pools gives an increasingly better average price than
using the Uniswap v3 USDC-WETH $.3\%$ fee tier pool alone. Specifically, we
route orders through the USDC-WETH $.3\%$, WETH-USDT $.3\%$, and USDC-USDT
$.01\%$ pools. This is the simplest example in which we can hope to achieve
improvements from routing, since two possible routes are available to the
seller: a direct route through the USDC-WETH pool; and an indirect route
that uses both the WETH-USDT pool and the USDC-USDT pool.

\section{Conclusion}
We constructed an efficient algorithm to solve the optimal routing problem. Our
algorithm parallelizes across markets and involves solving a series of optimal
arbitrage problems at each iteration. To facilitate efficient subproblem
solutions, we introduced an interface for swap markets, which includes aggregate
CFMMs.

We note that we implicitly assume that the trading sets are known
exactly when the routing problem is solved. This assumption, however, ignores 
the realities of trading on chain: unless our trades execute first in the next 
block, we are not guaranteed that the trading sets for each market are the same 
as those in the last block. Transactions before ours in the new block may have changed
prices (and reserves) of some of the markets we are routing through.
This observation naturally suggests \emph{robust routing} as
a natural direction for future research. Furthermore, efficient algorithms for 
routing with fixed transaction costs (\eg, gas costs) are another interesting 
direction for future work (see~\cite[\S 5]{angeris2022optimal} for the problem
formulation).

\section*{Acknowledgements}
We thank Francesco Iannelli and Jiahao Song for contributing to the
\texttt{CFMMRouter.jl} documentation and the Financial Cryptography 2023
reviewers for helpful comments.

\printbibliography

\appendix

\section{Closed form solutions}\label{app:closed-form}
Here, we cover some of the special cases where it is possible to analytically
write down the solutions to the arbitrage problems presented previously.

\paragraph{Geometric mean trading function.}
Some of the most popular swap markets, for example, Uniswap v2 and most
Balancer pools, which total over \$2B in reserves, are geometric mean
markets~\eqref{eq:geom} with $n=2$. This trading function can be written as
\[
    \phi(R) = R_1^wR_2^{1-w},
\]
where $0 < w < 1$ is a fixed parameter. This very common trading function
admits a closed-form solution to the arbitrage problem~\eqref{eq:arb}.
Using~\eqref{eq:func-form}, we can write
\[
    f_1(\delta_1) = R_2\left(1- \left(\frac{1}{1 + \gamma\delta_1/R_1}\right)^\eta\right)
\]
where $\eta = w/(1-w)$. (A similar equation holds for $f_2$.)
Using~\eqref{eq:swap-opt-cond} and~\eqref{eq:root}, and defining
\[
    \delta_1 = \frac{R_1}{\gamma}\left(\left(\eta\gamma\frac{\nu_2}{\nu_1}\frac{R_2}{R_1}\right)^{1/(\eta+1)} - 1\right),
\]
we have that $\delta_1^\star = \max\{\delta_1, 0\}$ is an optimal point
for~\eqref{eq:scalar-arb}. Note that when we take $w = 1/2$ then $\eta = 1$ and we recover the
optimal arbitrage for Uniswap given in~\cite[App.\ A]{angerisAnalysisUniswapMarkets2020}.

\paragraph{Bounded liquidity variation.} 
The bounded liquidity
variation~\eqref{eq:bounded-liquidity} of the product trading function
satisfies the definition of bounded liquidity given in~\S\ref{sec:bounded-liq}, whenever $\alpha,
\beta > 0$. 
We can write the forward exchange function
for the bounded liquidity product function~\eqref{eq:bounded-liquidity},
using~\eqref{eq:func-form}, as
\[
    f_1(\delta) = \min\left\{R_2, \frac{\gamma\delta(R_2 + \beta)}{R_1 +
        \gamma\delta + \alpha}\right\}
\]
The `min' here comes from the definition of a CFMM: it will not accept trades
which pay out more than the available reserves.
The maximum amount that a user can trade
with this market, which we will write as $\delta_1^-$, is when $f_1(\delta_1^-) = R_2$,
\ie,
\[
    \delta_1^- = \frac{1}{\gamma}\frac{R_2}{\beta}(R_1 + \alpha).
\]
(Note that this can also be derived by taking $f_1(\delta_1) = R_2$ in~\eqref{eq:func-form}
with the invariant~\eqref{eq:bounded-liquidity}.)
This means that
\[
    f_1^-(\delta_1^-) = \gamma\frac{\beta^2}{(R_1 + \alpha)(R_2 + \beta)},
\]
is the minimum supported price for asset 1. As before, a similar derivation
yields the case for asset 2.
Writing $k=(R_1+ \alpha)(R_2 + \beta)$, we see that we only need
to solve~\eqref{eq:scalar-arb} if the price $\nu_1/\nu_2$ is in the active
interval~\eqref{eq:interval},
\begin{equation}\label{eq:e-interval}
    \frac{\gamma \beta^2}{k} < \frac{\nu_1}{\nu_2} < \frac{k}{\gamma\alpha^2}.
\end{equation}
Otherwise, we know one of the two `boundary' solutions, $\delta_1^-$ or $\delta_2^-$, suffices.

\end{document}